\documentclass{paper}
\usepackage{amsmath,amsfonts}
\begin{document}
\title{The inhomogeneous invariance quantum group of q-deformed boson algebra with continuous
parameters}
\author{Azmi Ali Altintas$^1$, Metin Arik$^2$, Ali Serdar Arikan $^3$
\and $^1$%
{\small Okan University, Faculty of Engineering , Akf\i
rat, Istanbul, Turkey
 },\and
$^2${\small Bo\u{g}azi\c{c}i University, Department of Physics,
Bebek, Istanbul, Turkey} \and $^{3}${\small Sakarya University,
Department of Physics, Esentepe, Sakarya, Turkey}} \maketitle
e-mail:ali.altintas@okan.edu.tr \\
Keywords: Quantum Groups, q-deformed boson algebra.\\
\begin{abstract}
We present a q-deformed boson algebra using continuous momentum
parameters and investigate its inhomogeneous invariance quantum
group.
\end{abstract}
PACS:$02.20.Uw$

\section{Introduction}

Quantum group theory is related to integrable systems. To solve some
complicated problems classical groups are not sufficient and more
complex structures are used. These structures are quantum
groups and they are algebraic objects\cite{F,J}. There are many approaches to obtain quantum
groups. One way is called Drinfeld's approach\cite{D}. In that
approach by defining a deformation parameter on usual Lie algebra
one can construct a quantum group. The most familiar quantum group generated from a Lie algebra is $SU_q(2)$, here q is deformation parameter and it can take values between 0 and 1.
\begin{displaymath}
T=\left(
    \begin{array}{cc}
      a & b \\
      -qb^\ast & a^\ast \\
    \end{array}
  \right),
\end{displaymath} the elements satisfy
$ab=qba$, $ab^\ast =q b^\ast a$, $bb^\ast=b^\ast b$, $[a,a^\ast]=(q-1)bb^\ast$\\
 and the unitarity condition is: $TT^\ast=T^\ast T=\mathbf{1}$, also the q-determinant is equal to one:
$det_q(T)=aa^\ast +bb^\ast=a^\ast a+qb^\ast b=1$.

This procedure can be applied to some physical structure and one can get deformed physical systems. The most known example is deformed boson algebra.
By defining a deformation parameter on boson algebra one can write a deformed
boson algebra which is given by;
\begin{eqnarray}
a_ia_j^{\ast}-q_{ij}a_j^{\ast}a_i&=&\delta_{ij}\label{q1},\\
a_ia_j-a_ja_i&=&0\label{q2}.
\end{eqnarray}
where
\begin{equation}\label{q}
q_{ij}=(q-1)\delta_{ij}+1. \end{equation} Here i and j run from 1 to
n. Here $a\; \mbox{and}\; a^\ast$ are creation and annihilation
operators respectively.

 Another approach is realised by Woronowicz. The structures obtained are called matrix quantum groups. The
  idea of matrix quantum groups involves noncommutative comultiplication\cite{W}.Also Manin has shown that one can obtain quantum groups by
considering linear transformations on a quantum plane\cite{Manin}.
One can build up a matrix whose elements satisfy Hopf algebra
axioms. In this paper we use the term "invariance quantum group" to describe a Hopf algebra such that a particle algebra forms a right module
of the Hopf algebra.

Using Drinfeld's and Woronowicz's approach one can find inhomogeneous
invariance quantum groups of q deformed oscillators. This quantum
group is Bosonic Inhomogeneous Deformed General Linear Quantum
Group, $BIGL_{q_{ij},q_{ij}^{-1}}(2d)$ \cite{AAA}.

In quantum field theory the operators are not discrete instead they
are operators which have continuous parameters and the parameter is
usually momentum which is shown by $p$. In this manner boson algebra
can be written with continuous parameters, as;
\begin{eqnarray}
a(p)a^\ast(p')-a^\ast(p')a(p)&=&\delta(p-p'),\\
a(p)a(p')-a(p')a(p)&=&0.
\end{eqnarray}

 Now a question arises: if a q deformed oscillator is written using continuous
parameters what will be its inhomogeneous quantum invariance group?
In this study we will build a q oscillator algebra with continuous
parameters and find its inhomogeneous invariance quantum group.

\section{$BIGL_{g,g^{-1}}(2,\mathfrak{C})$}
Now we define a q oscillator algebra which the operators of the
algebra have continuous parameters. The bosonic creation and
annihilation operators are $a(p)$ and $a^{\ast}(p)$, here p is
thought as momentum and $p \in \mathbb{R}$. The q deformed boson
algebra which has continuous parameters can be written as
\begin{eqnarray}
a(p)a^{\ast}(p')-g(q;p,p')a^{\ast}(p')a(p)&=&\delta(p-p')\label{continuous1},\\
a(p)a(p')-a(p')a(p)&=&0\label{continuous2},
\end{eqnarray}
where
\begin{equation}
g(q;p,p')=\left\{\begin{array}{cc}
            q & p=p', \\
            1 & p \neq p'.
          \end{array}\right.
\end{equation}
It can be easily seen that the equations (\ref{continuous1}) and
(\ref{continuous2}) are compatible with the equations (\ref{q1}) and
(\ref{q2}).

Now we consider an inhomogeneous linear transformation on creation
and annihilation operators of the algebra. The transformed operators
can be written as;
\begin{eqnarray}
a(p)'&=&\int\alpha(p,k)\otimes a(k) dk+\int\beta(p,k)\otimes
a^{\ast}(k)
dk+f(p)\otimes 1,\\
a^{\ast}(p)'&=&\int\beta^\ast(p,k)\otimes a(k)
dk+\int\alpha^\ast(p,k)\otimes a^{\ast}(k) dk+ f^\ast(p)\otimes 1.
\end{eqnarray}
The parameters $\alpha(p,k)$ and $\beta(p,k)$ are not necessarily
commutative parameters and $f(p)$ is a noncommutative parameter. One
can write the transformation in matrix form.
\begin{equation}
M=\left(
\begin{array}{cc|c}
  \alpha  & \beta & f \\
  \beta^\ast & \alpha^\ast & f^\ast \\\hline
  0 & 0 & 1
\end{array}\right) =\left(
\begin{array}{c|c}
A&F\\\hline
0&1
\end{array}\right)
\end{equation}
Here A is homogeneous and F is inhomogeneous part of the
transformation matrix M. We want the transformed operators satisfy
the algebra which is defined in equations (\ref{continuous1}) and
(\ref{continuous2}). To leave the algebra invariant there should be
commutation relations between the elements of transformation matrix
M. The commutation relations
\begin{eqnarray}
\alpha(p,k)\alpha(p',l)&=&\alpha(p',l)\alpha(p,k),\\
\alpha(p,k)\alpha^\ast(p',l)&=&\frac{g(q;p,p')}{g(q;k,l)}\alpha^\ast(p',l)\alpha(p,k),\\
\alpha(p,k)\beta(p',l)&=&g(q;k,l)^{-1}\beta(p',l)\alpha(p,k),\\
\alpha(p,k)\beta^\ast(p',l)&=&g(q;p,p')\beta^{\ast}(p',l)\alpha(p,k),\\
\alpha(p,k)f(p')&=&f(p')\alpha(p,k),\\
\alpha(p,k)f^{\ast}(p')&=&g(q;p,p') f^{\ast}(p')\alpha(p,k),\\
\beta(p,k)\beta(p',l)&=&\beta(p',l)\beta(p,k),\\
\beta(p,k)\beta^{\ast}(p',l)&=&g(q;p,p')g(q;k,l)\beta^{\ast}(p',l)\beta(p,k),\\
\beta(p,k)f(p')&=&f(p')\beta(p,k),\\
\beta(p,k)f^{\ast}(p')&=&g(q;p,p') f^{\ast}(p')\beta(p,k),\\
f(p)f(p')-f(p')f(p)&=&\int\alpha(p,k)\beta(p',k)
dk-\int\alpha(p',k)\beta(p,k)dk,\\
f(p)f^{\ast}(p')-g(q;p,p')
f^{\ast}(p')f(p)&=&\delta(p-p')-\int\alpha(p,k)\alpha^{\ast}(p',k)dk+g(q;p,p')\int\beta^{\ast}(p',k)\beta(p,k)dk,\nonumber
\\
\end{eqnarray}
together with their hermitean conjugates.

In order to check that the transformation is a quantum group we
should find the coproduct, the counit and the coinverse of the
transformation. The coproduct of the transformation matrix M can be
found via tensor matrix multiplication.
\begin{equation}
\Delta(M)=M\dot{\otimes}M
\end{equation}
this gives;
\begin{eqnarray}
\Delta(\alpha(p,k))&=&\int\alpha(p,\eta)\otimes\alpha(\eta,k)d\eta+\int\beta(p,\eta)\otimes\beta^{\ast}(\eta,k)d\eta,\\
\Delta(\beta(p,k))&=&\int\alpha(p,\eta)\otimes\beta(\eta,k)d(\eta)+\int\beta(p,\eta)\otimes\alpha^{\ast}(\eta,k)d\eta,\\
\Delta(f(p))&=&\int\alpha(p,\eta)\otimes
f(\eta)d\eta+\int\beta(p,\eta)\otimes f^{\ast}(\eta) d(\eta)+f(p)\otimes 1,\\
\Delta(\alpha^{\ast}(p,k))&=&\int\alpha^{\ast}(p,\eta)\otimes\alpha^{\ast}(\eta,k)d\eta+\int\beta^{\ast}(p,\eta)\otimes\beta(\eta,k)d\eta,\\
\Delta(\beta^{\ast}(p,k))&=&\int\alpha^{\ast}(p,\eta)\otimes\beta^{\ast}(\eta,k)d\eta+\int\beta^{\ast}(p,\eta)\otimes\alpha(\eta,k)d\eta,\\
\Delta(f^{\ast}(p))&=&\int\alpha^{\ast}(p,\eta)\otimes
f^{\ast}(\eta)d\eta+\int\beta^{\ast}(p,\eta)\otimes f(\eta) d\eta+
f^{\ast}(p)\otimes 1,\\
\Delta(1)&=&1\otimes 1.
\end{eqnarray}
The counit and the coinverse are defined as;
\begin{eqnarray}
\varepsilon(M)&=&\mathcal{I},\\
S(M)&=&M^{-1}
\end{eqnarray}
In order to find the inverse of transformation matrix T we should
find the inverse of homogeneous part of transformation matrix.
\begin{equation}
T^{-1}=\left(
\begin{array}{cc}
A^{-1}&-A^{-1}F\\ 0&1
\end{array}\right),
\end{equation}
where
\begin{equation}
A=\left(\begin{array}{cc} \alpha(p,k)& \beta(p,k)\\
\beta^{\ast}(p,k)&\alpha^{\ast}(p,k)\end{array}\right)
\end{equation}
and
\begin{equation}
F=\left(\begin{array}{c} f(p)\\
f^{\ast}(p)\end{array}\right)
\end{equation}
 The commutation relations between the
elements of homogeneous part $A$ are;
\begin{eqnarray}
\alpha(p,k)\alpha(p',l)&=&\alpha(p',l)\alpha(p,k),\label{com1}\\
\alpha(p,k)\alpha^\ast(p',l)&=&\frac{g(q;p,p')}{g(q;k,l)}\alpha^\ast(p',l)\alpha(p,k),\\
\alpha(p,k)\beta(p',l)&=&g(q;k,l)^{-1}\beta(p',l)\alpha(p,k),\\
\alpha(p,k)\beta^\ast(p',l)&=&g(q;p,p')\beta^{\ast}(p',l)\alpha(p,k),\\
\beta(p,k)\beta(p',l)&=&\beta(p',l)\beta(p,k),\\
\beta(p,k)\beta^{\ast}(p',l)&=&g(q;p,p')g(q;k,l)\beta^{\ast}(p',l)\beta(p,k),\label{com6}
\end{eqnarray}
 We should remember that all $g(q;p,p')$'s are either $q$ or $1$.
 The equations \ref{com1} to \ref{com6} have two parameter deformed
 quantum group $GL_{p,q}$\cite{Sch} structure. One can write the inverse of
 $A$ in light of Schirrmacher's $GL_{p,q}$. In our case
 $p=g(q;k,l)^{-1}$ and $q=g(q;p,p')^{-1}$.

 Since the elements of matrix T have coproduct, counit and coinverse
 we can say that this transformation is a quantum group. This group
 called as $BIGL_{g,g^{-1}}(2,\mathfrak{C})$ Bosonic Inhomogeneous Deformed General Linear
 Quantum Group with continuous parameter.
 \section{Conclusion}
 In this study we have shown that there is an inhomogeneous
 invariance quantum group for deformed boson algebra with continuous
 parameters. This quantum group is $BIGL_{g,g^{-1}}(2)$. It is easy
 to see that the inhomogeneous part has bosonic structure.

 The homogeneous part of $BIGL_{g,g^{-1}}(2,\mathfrak{C})$ has also a quantum
 group structure. This is $GL_{g,g^{-1}}(2)$. If one can choose $g=1$
 the quantum group becomes $BISp(2,\mathfrak{C})$, Bosonic Inhomogeneous
 Symplectic Quantum Group with continuous parameter\cite{Ali}.

 In field theory a field $\varphi(x)$ can be written using creation
 and annihilation operators. Using the deformed creation and
 annihilation operators $a(p)$ and $a^{\ast}(p)$ one can construct a
 deformed field theory. Also using transformed creation and
 annihilation operators one can find a quantum group for deformed
 quantum field theory. The theory of quantum groups can be helpful in the generalization of quantization and a more consistent approach to
 interacting field theories.


\begin{thebibliography}{99}
\bibitem{F}
L.D. Faddeev, N.Y. Reshetikhin, L.A. Takhtajan, Leningrad Math J.
\textbf{1}, 193 (1990)
\bibitem{J}
M. Jimbo, Lett. Math. Phys. \textbf{11}, 247 (1986)
\bibitem{D}
V.G. Drinfeld, Proceedings International Congr. Math., Berkeley
\textbf{1}, 798 (1986)
\bibitem{W}
S.L. Woronowicz, Commun. Math. Phys. \textbf{111}, 613 (1987)
\bibitem{Manin}
Yu I. Manin, Quantum Groups and Noncommutative Geometry, Centre de
Reserches Mathematiques, Montreal, (1988)
\bibitem{AAA}
A.A. Altintas, M. Arik, A.S. Arikan, CEJP. 8(1), 131, (2010)
\bibitem{Sch}
A. Schirrmacher, J. Wess, B. Zumino, Z.Phys. C-Particles and Fields
49, 317 (1991).
\bibitem{Ali}
A.A. Altintas, M. Arik, Int. J. Theor. Phys. 48,1534, (2009).
\end{thebibliography}
\end{document}